\documentclass[12pt]{amsart}

\usepackage{amsmath}
\usepackage{amssymb}
\usepackage{url}
\usepackage{amscd}
\usepackage{mathrsfs}
\usepackage[dvips]{graphicx}
\usepackage{indentfirst}
\usepackage{setspace}

\theoremstyle{plain}
\newtheorem{thm}{Theorem}

\newtheorem{lemma}{Lemma}

\newcommand{\Z}{{\mathbb Z}}

\newcommand{\gf}{{\mathrm{GF}}}

\newcommand{\cc}{{\mathbf c}}

\newcommand{\tr}{{\mathrm{Tr}}}

\newcommand{\F}{{\mathcal F}}
\newcommand{\A}{{\mathcal A}}

\newcommand{\Ta}{{\mathcal T}_1}
\newcommand{\Tb}{{\mathcal T}_2}

\newcommand{\FF}{{\mathbb{F}}}

\newcommand{\C}{{\mathcal C}}

\renewcommand{\vec}[1]{\underline{#1}}

\def\({\left(}
\def\){\right)}

\begin{document}

\title{The weight distributions of a class of cyclic codes II}

\author[M. Xiong]{\sc Maosheng Xiong}

\address{Maosheng Xiong: Department of Mathematics,
Hong Kong University of Science and Technology,
Clear Water Bay, Kowloon, Hong Kong
}
\email{mamsxiong@ust.hk}

\thanks{The author was supported by the Research Grants Council of Hong Kong under Project Nos. RGC606211 and DAG11SC02.}


\begin{abstract}
Recently, the weight distributions of the duals of the cyclic codes with two zeros have been obtained for several cases in \cite{DL1,DL2,WT,X}. In this paper we use the method developed in \cite{X} to solve one more special case. We make extensive use of standard tools in number theory such as characters of finite fields, the Gauss sums and the Jacobi sums. The problem of finding the weight distribution is transformed into a problem of evaluating certain character sums over finite fields, which turns out to be associated with counting the number of points on some elliptic curves over finite fields. We also treat the special case that the characteristic of the finite field is 2. 
\end{abstract}

\maketitle

\thispagestyle{empty}

\maketitle

\thispagestyle{empty}

\section{Introduction}
Denote by $\gf(q)$ the finite field of order $q$, where $q=p^s$, $s$ is a positive integer and $p$ is a prime number. An $[n,k,d]$-linear code $\C$ over $\gf(q)$ is a $k$-dimensional subspace of $\gf(q)^n$ with minimum distance $d$. If, in addition, $\C$ satisfies the condition that $(c_{n-1},c_0,c_1,\ldots,c_{n-2}) \in \C$ whenever $(c_0,c_1,\ldots,c_{n-2},c_{n-1}) \in \C$, then $\C$ is a cyclic code. Let $A_i$ denote the number of codewords with Hamming weight $i$ in $\C$. The weight enumerator of $\C$ is defined by
\[1+A_1x+A_2x^2+\cdots+A_nx^n.\]
The sequence $(1,A_1,\ldots,A_n)$ is called the weight distribution of $\C$. In coding theory it is often desirable to know the weight distribution of a code because it contains a lot of important information, for example, it can be used to estimate the error correcting capability and the probability of error detection and correction with respect to some algorithms. This is quite useful in practice. Many important families of cyclic codes have been studied extensively in the literature, so have their various properties. However the weight distributions are difficult to obtain in general and they are known only for a few special families.

Given a positive integer $m$, let $r=q^m$, and $\alpha$ be a generator of the multiplicative group $\gf(r)^*:=\gf(r)-\{0\}$. Let $h$ be a positive factor of $q-1$ and $1<e$ be an integer such that $e|\gcd\left({q-1},hm\right)$. We define
\begin{eqnarray} \label{1:para} g=\alpha^{(q-1)/h}, n=\frac{h(r-1)}{q-1}, \beta=\alpha^{(r-1)/e}, N=\gcd\left(\frac{r-1}{q-1},\frac{e(q-1)}{h}\right).\end{eqnarray}
It is clear that the order of $g$ is $n$ and $(g \beta)^n=1$. Refining an argument from \cite{DL1}, we will prove later that the minimal polynomials of $g^{-1}$ and $(\beta g)^{-1}$ are distinct over $\gf(q)$, hence their product is a factor of $x^n-1$, except when $q=3,h=1,e=m=2$. We remark that the conditions here are slightly more general than those in \cite{DL1,DL2,X} (and in several other references), which require that $e|h$. This consideration is inspired by an anonymous referee and actually provides more flexibility.

Define the cyclic code over $\gf(q)$ by
\begin{eqnarray} \label{1:code} \C_{(q,m,h,e)}=\left\{\cc_{(a,b)}: a,b \in \gf(r)\right\},\end{eqnarray}
where the codeword $\cc_{(a,b)}$ is give by
\begin{eqnarray} \label{1:codedef} \cc_{(a,b)}:=\left(\tr\left(ag^i+b(\beta g)^i\right)\right)_{i=0}^{n-1}.\end{eqnarray}
Here for simplicity $\tr$ is the trace function from $\gf(r)$ to $\gf(q)$.

The code $\C_{(q,m,h,e)}$ has been an interesting subject of study for a long time. For example, when $h=q-1$, the code $\C_{(q,m,h,e)}$ is the dual of a primitive cyclic linear code with two zeros; such codes have been studied extensively (see for example \cite{BM, CCD, CCZ,C,LF1, LTW, Mc,MR,S,YCD}). In general the dimension of $\C_{(q,m,h,e)}$ is $2m$, but determining the weight distribution is very difficult. However, in certain special cases the weight distribution is known. We summarize these cases below.
\begin{itemize}
\item[1)] $e>1$ and $N=1$ (\cite{DL1});

\item[2)] $e=2$ and $N=2$ (\cite{DL1});

\item[3)] $e=2$ and $N=3$ (\cite{DL2});

\item[4)] $e=2$ and $p^j+1 \equiv 0 \pmod{N}$, where $j$ is a positive integer (\cite{DL2});

\item[5)] $e=3$ and $N=2$ (\cite{WT});

\item[6)] $e=4$ and $N=2$ (\cite{X}).

\end{itemize}

In this paper we compute the weight distribution for one more case $e=N=3$. As it turns out, if $p \equiv 1 \pmod{3}$, the number of distinct nonzero weights in the codes is 12 or 13, the shortest code in the family has length $\frac{p^3-1}{p-1}$ over $\gf(p)$. On the other hand, if $p \equiv 2 \pmod{3}$, then the number of distinct nonzero weights is $5$ or $6$, and the shortest code in the family has length $\frac{p^6-1}{p^2-1}$ over $\gf(p^2)$. The dimension is always $2m$ where $3|m$, so the smallest dimension is $6$. We have computed several examples for relatively small parameters by Magma, and thanks to the referee's suggestions, we also compare them with the best existing codes from Markus Grassl's table (http://www.codetables.de/). It seems the codes constructed in this way fall short of such comparison.

To describe the results, for the sake of clear presentation, we use the ``modified'' weight $\lambda(a,b)$, instead of the usual Hamming weight $w(\cc_{(a,b)})$ for a codeword $\cc_{(a,b)}$. The relation between them is given by the formula
\begin{eqnarray} \label{1:modweight}
w(\cc_{(a,b)})=\frac{h(r-1)}{q}-\lambda(a,b).
\end{eqnarray}
The case that $p \equiv 2 \pmod{3}$ is easy to describe.
\begin{thm} \label{thm1} Let $\C_{(q,m,h,e)}$ be the cyclic code defined by (\ref{1:code}) and (\ref{1:codedef}), and the parameters are given by (\ref{1:para}) where $q=p^s$. Assume that $e=N=3$ and $p \equiv 2 \pmod{3}$ (including the case $p=2$).
\begin{itemize}
\item[(1).] If $3|\frac{q-1}{h}$, the modified weight distribution of $\C_{(q,m,h,e)}$ is given by Table \ref{1:t1}.

\item[(2).] If $3 \nmid \frac{q-1}{h}$, the modified weight distribution of $\C_{(q,m,h,e)}$ is given by Table \ref{1:t2}.

\end{itemize}
\end{thm}

\begin{table}[ht] 
\caption{The case $e=N=3$, $p \equiv 2 \pmod{3}$ and $3|\frac{q-1}{h}$}
\centering
\begin{tabular}{|c|| c|}
\hline
\hline
Weight $\lambda(a,b)$ & Frequency  \\
[0.5ex]
\hline
$-\frac{h}{q}\left\{2(-1)^{ms/2}\sqrt{r}+1\right\}$ & $\frac{r-1}{27}\bigl\{r-8-2(-1)^{ms/2}\sqrt{r}\bigr\}$ \\
\hline
$\frac{h}{q}\left\{(-1)^{ms/2}\sqrt{r}-1\right\}$ & $\frac{2(r-1)}{27}\bigl\{4r-14+(-1)^{ms/2}\sqrt{r}\bigr\}$ \\
\hline
$-\frac{h}{q}\left\{(-1)^{ms/2}\sqrt{r}+1\right\}$ & $\frac{2(r-1)}{9}\bigl\{r-2+(-1)^{ms/2}\sqrt{r}\bigr\}$ \\
\hline
$-\frac{h}{q}$ & $\frac{2(r-1)}{9}\bigl\{2r-1-(-1)^{ms/2}\sqrt{r}\bigr\}$ \\
\hline
$\frac{h}{3q}\left\{r-4(-1)^{ms/2}\sqrt{r}-3\right\}$ & $r-1$ \\
\hline
$\frac{h}{3q}\left\{r+2(-1)^{ms/2}\sqrt{r}-3\right\}$ & $2(r-1)$ \\
\hline
$\frac{h(r-1)}{q}$ & $1$ \\
\hline
\hline
\end{tabular}
\label{1:t1}
\end{table}

\begin{table}[ht] 
\caption{The case $e=N=3$, $p \equiv 2 \pmod{3}$, and $3 \nmid \frac{q-1}{h}$}
\centering
\begin{tabular}{|c| c|}
\hline
\hline
Weight $\lambda(a,b)$ & Frequency  \\
[0.5ex]
\hline
$-\frac{h}{q}\left\{2(-1)^{ms/2}\sqrt{r}+1\right\}$ & $\frac{r-1}{27}\bigl\{r+1-2(-1)^{ms/2}\sqrt{r}\bigr\}$ \\
\hline
$\frac{h}{q}\left\{(-1)^{ms/2}\sqrt{r}-1\right\}$ & $\frac{2(r-1)}{27}\bigl\{4r-5+(-1)^{ms/2}\sqrt{r}\bigr\}$ \\
\hline
$-\frac{h}{q}\left\{(-1)^{ms/2}\sqrt{r}+1\right\}$ & $\frac{2(r-1)}{9}\bigl\{r-2+(-1)^{ms/2}\sqrt{r}\bigr\}$ \\
\hline
$-\frac{h}{q}$ & $\frac{r-1}{9}\bigl\{4r-11-2(-1)^{ms/2}\sqrt{r}\bigr\}$ \\
\hline
$\frac{h}{3q}\left\{r+2(-1)^{ms/2}\sqrt{r}-3\right\}$ & $r-1$ \\
\hline
$\frac{h}{3q}\left\{r-(-1)^{ms/2}\sqrt{r}-3\right\}$ & $2(r-1)$ \\
\hline
$\frac{h(r-1)}{q}$ & $1$ \\
\hline
\hline
\end{tabular}
\label{1:t2}
\end{table}

\noindent {\bf Example 1.} Let $p=2,s=2,q=4,m=3,r=64,h=1, e=N=3$. Letting $\alpha$ be a generator of $\gf(64)$ from Magma, which uses the irreducible polynomial $x^6+x^4+x^3+x+1$, we can construct the code explicitly and the weight distribution of the cyclic code $\C_{(q,m,h,e)}$ is given by
\[1+63x^{8}+294x^{12}+756x^{14}+1890x^{16}+1092x^{18}. \]
This is confirmed by computing Table \ref{1:t1}, since $3 | \frac{q-1}{h}=3$. Notice that there are only six weights because two of the weights in Table \ref{1:t1} are the same, namely,
\[-\frac{h\left(2(-1)^{ms/2}\sqrt{r}+1\right)}{q}=\frac{h\left(r+2(-1)^{ms/2}\sqrt{r}-3\right)}{3q}.\]
Actually this happens if and only if $r=2^6=64$. In other cases, there are always seven distinct weights.

This is a $[21,6,8]$-cyclic code over $\gf(4)$. Grassl's table shows that there is a $[21,6,12]$ code over $\gf(4)$, and the best possible minimum distance is 12. $\quad \square$

\noindent {\bf Example 2.} Let $p=2,s=2,q=4,m=3,r=64,e=h=3,N=3$. Letting $\alpha$ be a generator of $\gf(64)$ from Magma, which uses the irreducible polynomial $x^6+x^4+x^3+x+1$, we can construct the code explicitly and the weight distribution of the cyclic code $\C_{(q,m,h,e)}$ is given by 
\[1+126x^{30}+252x^{36}+756x^{42}+1827x^{48}+1134x^{54}. \]
This is confirmed by computing Table \ref{1:t2}, since $3 \nmid \frac{q-1}{h}=1$. Same as Example 1, there are only six weights because two of the weights in Table \ref{1:t2} are the same.

This is a $[63,6,30]$-cyclic code over $\gf(4)$. Grassl's table shows that there is a $[63,6,44]$ code over $\gf(4)$, and the best possible minimum distance is 44. $\quad \square$

The case that $p \equiv 1 \pmod{3}$ can be described but the results are a little more complicated, because they rely on a subtle choice of cubic characters of $\gf(p)$ and $\gf(r)$, which need to be made explicit. We list the results first and then explain how to compute them later.

\begin{thm} \label{thm2} Let $\C_{(q,m,h,e)}$ be the cyclic code defined by (\ref{1:code}) and (\ref{1:codedef}), and the parameters are given by (\ref{1:para}) where $q=p^s$. Assume that $e=N=3$ and $p \equiv 1 \pmod{3}$.
\begin{itemize}
\item[(1).] If $3|\frac{q-1}{h}$, the modified weight distribution of $\C_{(q,m,h,e)}$ is given by Table \ref{1:t3}.

\item[(2).] If $3 \nmid \frac{q-1}{h}$, the modified weight distribution of $\C_{(q,m,h,e)}$ is given by Table \ref{1:t4}.
\end{itemize}
The symbols in Tables \ref{1:t3}-\ref{1:t4} are as follows: let \[\omega:=\frac{-1+\sqrt{-3}}{2}.\]
Then there is a unique algebraic integer $\pi \in \Z[\omega]$ such that $\pi \bar{\pi}=p$ and $\pi \equiv -1 \pmod{3}$ ($\bar{\pi}$ is the complex conjugate of $\pi$), and we let $\rho$ be the cubic character of $\gf(r)$ arising from $\left(\frac{\cdot}{\pi}\right)_3$, the standard cubic residue symbol of the ring $\Z[\omega]/\pi\Z[\omega]$, which is isomorphic to $\gf(p)$, and the three Gaussian periods $\eta_1^{(3,r)},\eta_{\alpha}^{(3,r)}$ and $\eta_{\alpha^2}^{(3,r)}$ are given by
\begin{eqnarray*} \eta_1^{(3,r)}&=&\frac{(-1)^{sm+1}r^{1/3}\left(\pi^{sm/3}+\bar{\pi}^{sm/3}
\right)-1}{3},\\
\eta_{\alpha}^{(3,r)}&=&\frac{(-1)^{sm+1}r^{1/3}\left(\rho^2(\alpha)\pi^{sm/3}+\rho(\alpha)\bar{\pi}^{sm/3}
\right)-1}{3},\\
\eta_{\alpha^2}^{(3,r)}&=&\frac{(-1)^{sm+1}r^{1/3}\left(\rho(\alpha)\pi^{sm/3}+\rho^2(\alpha)\bar{\pi}^{sm/3}
\right)-1}{3}.
\end{eqnarray*}
Moreover, the relation between the modified weight $\lambda(a,b)$ and the Hamming weight $w(\cc_{(a,b)})$ is given by (\ref{1:modweight}).
\end{thm}

For computational purposes, we can choose the value of $\pi$ explicitly. By the unique factorization property of the ring $\Z[\omega]$, the prime $p$ can be represented as $p=a^2-ab+b^2$ for some integers $a,b$. If, in addition, we require that $a \equiv 2 \pmod{3}$, $b \equiv 0 \pmod{3}$ and $b>0$, then such integers $a,b$ exist and are unique. We can choose $\pi=a+b \omega$. Hence $\bar{\pi}=a+b \omega^2=(a-b)-b\omega$.

We still need to determine the value of $\rho(\alpha)$, which is either $\omega$ or $\omega^2$. This can be done by using the definition of $\left(\frac{\cdot}{\pi}\right)_3$ and the explicit identification of $\Z[\omega]/\pi\Z[\omega]$ with $\gf(p)$ (see \cite[Chapter 9]{IR}), and we can describe the algorithm as follows: since $p \nmid b$, there is an integer $b'$ such that $bb' \equiv 1 \pmod{p}$. Let $N_{r/p}:\gf(r) \to \gf(p)$ be the norm map, we know that $N_{r/p}(\alpha)=\alpha^{(r-1)/(p-1)} \in \gf(p)$, and $N_{r/p}(\alpha)$ can be naturally identified by an integer modulo $p$. It can be shown that either $N_{r/p}(\alpha)^{(p-1)/3} \equiv -b'a \pmod{p}$ or $N_{r/p}(\alpha)^{(p-1)/3} \equiv -1+b'a \pmod{p}$. Then the value of $\rho(\alpha)$ is given by
\[\rho(\alpha)=\left\{\begin{array}{lll}
\omega&:& \mbox{ if } N_{r/p}(\alpha)^{(p-1)/3} \equiv -b'a \pmod{p},\\
\omega^2&:& \mbox{ if } N_{r/p}(\alpha)^{(p-1)/3} \equiv -1+b'a \pmod{p}.
\end{array}\right.
\]

\noindent {\bf Example 3.} Let $p=q=7,s=1,m=3,h=1,e=N=3$. Letting $\alpha$ be a generator of $\gf(7^3)$ from Magma, which uses the irreducible polynomial $x^3+6x^2+4$, we can construct the code explicitly and find that the weight distribution of the cyclic code $\C_{(q,m,h,e)}$ is
\begin{eqnarray*} 1+342x^{30}+342x^{32}+342x^{36}+3990x^{45}+
14364x^{46}+12312x^{47}\\
+16302x^{48}+24624x^{49}
+14364x^{50}+14364x^{51}+12312x^{52}+3990x^{54}. \end{eqnarray*}
This is confirmed by computing Table \ref{1:t3}, since $p \equiv 1 \pmod{3}$ and $3 | \frac{q-1}{h}=6$. Here $7=2^2-2*3+3^2$, so $a=2,b=3,\pi=2+3 \omega$. We choose $b'=-2$ so that $bb' \equiv 1 \pmod{7}$. We find from Magma that $N_{r/p}(\alpha)=3$. Hence
$3^{(p-1)/3}=3^2 \equiv 2 \pmod{7}$ and $-1+a b'=-1+2*(-2)=-5 \equiv 2 \pmod{7}$, so $\rho(\alpha)=\omega^2$. The Gaussian periods can be computed as $\eta_1^{(3,r)}=2, \eta_{\alpha}^{(3,r)}=-12,\eta_{\alpha^2}^{(3,r)}=9$. Now the weight distribution can be obtained from Table \ref{1:t3}. There are only 13 weights because two of the weights in Table \ref{1:t3} are the same, namely
\[\frac{3h\eta_1^{(3,r)}}{q}=\frac{h(2\eta_{\alpha^2}^{(3,r)}+\eta_{\alpha}^{(3,r)})}{q}=\frac{18}{7}. \]

This is a $[57,6,30]$-cyclic code over $\gf(7)$. Grassl's table shows that there is a $[57,6,42]$ code over $\gf(7)$, and the best possible minimum distance can not be larger than 45. $\quad \square$

\noindent {\bf Example 4.} Let $p=q=7,s=1,m=3,e=h=3,N=3$. Letting $\alpha$ be a generator of $\gf(7^3)$ from Magma, which uses the irreducible polynomial $x^3+6x^2+4$, we can construct the code explicitly and find that the weight distribution of the cyclic code $\C_{(q,m,h,e)}$ is
\begin{eqnarray*} 1+342x^{93}+342x^{99}+342x^{102}+4104x^{135}+
14364x^{138}+12312x^{141}\\
+16416x^{144}+24282x^{147}
+14364x^{150}+14364x^{153}+12312x^{156}+4104x^{162}. \end{eqnarray*}
This is confirmed by computing Table \ref{1:t4}, since $p \equiv 1 \pmod{3}$ and $3 \nmid \frac{q-1}{h}=2$. As in Example 3, we find $\rho(\alpha)=\omega^2$, $\eta_1^{(3,r)}=2, \eta_{\alpha}^{(3,r)}=-12$ and $\eta_{\alpha^2}^{(3,r)}=9$. Now the weight distribution can be obtained from Table \ref{1:t4}. Same as Example 3, there are only 13 weights because two of the weights in Table \ref{1:t4} are the same. 

This is a $[171,6,93]$-cyclic code over $\gf(7)$. The length 171 of the code is too large for comparison with Grassl's table. $\quad \square$


\begin{table}[ht] 
\caption{The case $e=N=3$, $p \equiv 1 \pmod{3}$, and $3|\frac{q-1}{h}$}
\centering
\begin{tabular}{|c| c|}
\hline
\hline
Weight $\lambda(a,b)$ & Frequency  \\
[0.5ex]
\hline
$\frac{3h}{q} \eta_1^{(3,r)}$ & $\frac{r-1}{27}\bigl\{r-8-(-1)^{ms}\left(\pi^{ms}+\bar{\pi}^{ms}\right)\bigr\}$ \\
\hline
$\frac{3h}{q} \eta_{\alpha}^{(3,r)}$ & $\frac{r-1}{27}\bigl\{r-8-(-1)^{ms}\left(\pi^{ms}+\bar{\pi}^{ms}\right)\bigr\}$ \\
\hline
$\frac{3h}{q} \eta_{\alpha^2}^{(3,r)}$ & $\frac{r-1}{27}\bigl\{r-8-(-1)^{ms}\left(\pi^{ms}+\bar{\pi}^{ms}\right)\bigr\}$ \\
\hline
$\frac{h}{q} \left\{2\eta_1^{(3,r)}+\eta_{\alpha}^{(3,r)}\right\}$ & $\frac{r-1}{9}\bigl\{r-2-(-1)^{ms}\left(\rho^2(\alpha)\pi^{ms}+\rho(\alpha)\bar{\pi}^{ms}\right)\bigr\}$ \\
\hline
$\frac{h}{q} \left\{2\eta_1^{(3,r)}+\eta_{\alpha^2}^{(3,r)}\right\}$ & $\frac{r-1}{9}\bigl\{r-2-(-1)^{ms}\left(\rho(\alpha)\pi^{ms}+\rho^2(\alpha)\bar{\pi}^{ms}\right)\bigr\}$ \\
\hline
$\frac{h}{q} \left\{2\eta_{\alpha}^{(3,r)}+\eta_{1}^{(3,r)}\right\}$ & $\frac{r-1}{9}\bigl\{r-2-(-1)^{ms}\left(\rho(\alpha)\pi^{ms}+\rho^2(\alpha)\bar{\pi}^{ms}\right)\bigr\}$ \\
\hline
$\frac{h}{q} \left\{2\eta_{\alpha}^{(3,r)}+\eta_{\alpha^2}^{(3,r)}\right\}$ & $\frac{r-1}{9}\bigl\{r-2-(-1)^{ms}\left(\rho^2(\alpha)\pi^{ms}+\rho(\alpha)\bar{\pi}^{ms}\right)\bigr\}$ \\
\hline
$\frac{h}{q} \left\{2\eta_{\alpha^2}^{(3,r)}+\eta_{1}^{(3,r)}\right\}$ & $\frac{r-1}{9}\bigl\{r-2-(-1)^{ms}\left(\rho^2(\alpha)\pi^{ms}+\rho(\alpha)\bar{\pi}^{ms}\right)\bigr\}$ \\
\hline
$\frac{h}{q} \left\{2\eta_{\alpha^2}^{(3,r)}+\eta_{\alpha}^{(3,r)}\right\}$ & $\frac{r-1}{9}\bigl\{r-2-(-1)^{ms}\left(\rho(\alpha)\pi^{ms}+\rho^2(\alpha)\bar{\pi}^{ms}\right)\bigr\}$ \\
\hline
$-\frac{h}{q}$ & $\frac{2(r-1)}{9}\bigl\{r+1-(-1)^{ms}\left(\pi^{ms}+\bar{\pi}^{ms}\right)\bigr\}$ \\
\hline
$\frac{h}{3q}\left\{r-1+6 \eta_{1}^{(3,r)}\right\}$ & $r-1$ \\
\hline
$\frac{h}{3q}\left\{r-1+6 \eta_{\alpha}^{(3,r)}\right\}$ & $r-1$ \\
\hline
$\frac{h}{3q}\left\{r-1+6 \eta_{\alpha^2}^{(3,r)}\right\}$ & $r-1$ \\
\hline
$\frac{h(r-1)}{q}$ & $1$ \\
\hline
\hline
\end{tabular}
\label{1:t3}
\end{table}

\begin{table}[ht] 
\caption{The case $e=N=3$, $p \equiv 1 \pmod{3}$, and $3 \nmid \frac{q-1}{h}$}
\centering
\begin{tabular}{|c| c|}
\hline
\hline
Weight $\lambda(a,b)$ & Frequency  \\
[0.5ex]
\hline
$\frac{3h}{q} \eta_1^{(3,r)}$ & $\frac{r-1}{27}\bigl\{r+1-(-1)^{ms}\left(\pi^{ms}+\bar{\pi}^{ms}\right)\bigr\}$ \\
\hline
$\frac{3h}{q} \eta_{\alpha}^{(3,r)}$ & $\frac{r-1}{27}\bigl\{r+1-(-1)^{ms}\left(\pi^{ms}+\bar{\pi}^{ms}\right)\bigr\}$ \\
\hline
$\frac{3h}{q} \eta_{\alpha^2}^{(3,r)}$ & $\frac{r-1}{27}\bigl\{r+1-(-1)^{ms}\left(\pi^{ms}+\bar{\pi}^{ms}\right)\bigr\}$ \\
\hline
$\frac{h}{q} \left\{2\eta_1^{(3,r)}+\eta_{\alpha}^{(3,r)}\right\}$ & $\frac{r-1}{9}\bigl\{r-2-(-1)^{ms}\left(\rho^2(\alpha)\pi^{ms}+\rho(\alpha)\bar{\pi}^{ms}\right)\bigr\}$ \\
\hline
$\frac{h}{q} \left\{2\eta_1^{(3,r)}+\eta_{\alpha^2}^{(3,r)}\right\}$ & $\frac{r-1}{9}\bigl\{r-2-(-1)^{ms}\left(\rho(\alpha)\pi^{ms}+\rho^2(\alpha)\bar{\pi}^{ms}\right)\bigr\}$ \\
\hline
$\frac{h}{q} \left\{2\eta_{\alpha}^{(3,r)}+\eta_{1}^{(3,r)}\right\}$ & $\frac{r-1}{9}\bigl\{r-2-(-1)^{ms}\left(\rho(\alpha)\pi^{ms}+\rho^2(\alpha)\bar{\pi}^{ms}\right)\bigr\}$ \\
\hline
$\frac{h}{q} \left\{2\eta_{\alpha}^{(3,r)}+\eta_{\alpha^2}^{(3,r)}\right\}$ & $\frac{r-1}{9}\bigl\{r-2-(-1)^{ms}\left(\rho^2(\alpha)\pi^{ms}+\rho(\alpha)\bar{\pi}^{ms}\right)\bigr\}$ \\
\hline
$\frac{h}{q} \left\{2\eta_{\alpha^2}^{(3,r)}+\eta_{1}^{(3,r)}\right\}$ & $\frac{r-1}{9}\bigl\{r-2-(-1)^{ms}\left(\rho^2(\alpha)\pi^{ms}+\rho(\alpha)\bar{\pi}^{ms}\right)\bigr\}$ \\
\hline
$\frac{h}{q} \left\{2\eta_{\alpha^2}^{(3,r)}+\eta_{\alpha}^{(3,r)}\right\}$ & $\frac{r-1}{9}\bigl\{r-2-(-1)^{ms}\left(\rho(\alpha)\pi^{ms}+\rho^2(\alpha)\bar{\pi}^{ms}\right)\bigr\}$ \\
\hline
$-\frac{h}{q}$ & $\frac{r-1}{9}\bigl\{2r-7-2(-1)^{ms}\left(\pi^{ms}+\bar{\pi}^{ms}\right)\bigr\}$ \\
\hline
$\frac{h}{q}\left(\frac{r-1}{3}+ \eta_{1}^{(3,r)}+\eta_{\alpha}^{(3,r)}\right)$ & $r-1$ \\
\hline
$\frac{h}{q}\left(\frac{r-1}{3}+ \eta_{1}^{(3,r)}+\eta_{\alpha^2}^{(3,r)}\right)$ & $r-1$ \\
\hline
$\frac{h}{q}\left(\frac{r-1}{3}+\eta_{\alpha}^{(3,r)}+ \eta_{\alpha^2}^{(3,r)}\right)$ & $r-1$ \\
\hline
$\frac{h(r-1)}{q}$ & $1$ \\
\hline
\hline
\end{tabular}
\label{1:t4}
\end{table}

The ideas of the proofs of Theorems 1 and 2 are similar to those of \cite{X}, that is, first we use orthogonal properties of characters to transform the problem of finding the weight distribution into a problem of evaluating certain character sums over finite fields, and then we group character sums accordingly and relate them to counting the number of points on some curves over finite fields. For $e=N=3$, the curves turn out to be elliptic curves, on which the number of points can be computed explicitly by using standard techniques involving Gauss sums and Jacobi sums. While the methods are similar, the problem in this paper is more complicated, because firstly $N=3$, so there are three Gaussian periods instead of two; the difference is significant. Secondly, the prime $p$ could be $p=2$, $p \equiv 1 \pmod{3}$ and $p \equiv 2 \pmod{3}$, all of which need to be taken care of; moreover, whether or not 3 divides $\frac{q-1}{h}$ also has an effect. There are simply quite a few cases to consider and a lot of computation is involved. It turns out that the case $p=2$ can be included into the case $p \equiv 2 \pmod{3}$.

The paper is organized as follows: in Section 2 we recall the result we obtained in \cite{X} and apply it to the case $e=N=3$; the cases $p=2$ and $p$ odd need to be taken care of separately; Section 3 is devoted to the proof of Theorem 2; in Section 4, we argue that $p=2$ can be included into the case $p \equiv 2 \pmod{3}$, and we prove Theorem 1. We find the papers \cite{DL1,DL2,WT} quite inspiring and very well-written, which we use as general references and starting points of this paper. Interested readers may refer to them for some preliminary background and other information related to this subject.

\noindent {\bf Acknowledgment} The author is grateful to the anonymous referees for many valuable suggestions which help improve the quality of the paper substantially, and to Professor Cunsheng Ding for bringing this problem to his attention.

\section{Preliminaries}

We first prove that, given parameters $p,q,r,m,s,\beta,g, \ldots$ etc with $r=q^m$ in (\ref{1:para}) and conditions $h|(q-1)$ and $1<e|\gcd(q-1,hm)$, the order of $g$ is $n$, $(g \beta)^n=1$ and the minimal polynomials of $g^{-1}$ and $(\beta g)^{-1}$ are distinct over $\gf(q)$ except when $q=3,h=1,e=m=2$.

It is clear that the order of $g$ is $n$. Since $q \equiv 1 \pmod{e}$ and
\[\frac{h(r-1)}{q-1}=h(q^{m-1}+\cdots+q+1) \equiv hm \equiv 0 \pmod{e},\] 
we also have $(g \beta)^n=1$. 

Now suppose that the minimal polynomials of $g^{-1}$ and $(\beta g)^{-1}$ are the same over $\gf(q)$, then there is an integer $a$, $0 \le a \le m-1$ such that
\[\frac{q-1}{h}(q^a-1) \equiv \frac{q^m-1}{e} \pmod{q^m-1}. \]
Clearly $a \ne 0$ as $e>1$, so $a \ge 1$ and $m \ge 2$. Since 
\[0<\frac{q-1}{h}(q^a-1), \frac{q^m-1}{e} <q^m-1,\] 
we have the equality
\begin{eqnarray} \label{1:eqpa} \frac{q-1}{h}(q^a-1) = \frac{q^m-1}{e}. \end{eqnarray}
If $a \le m-2$, then
\[\frac{(q-1)e}{h}(q^a-1) \le q^2(q^{m-2}-1)=q^m-q^2<q^m-1,\]
so we must have $a=m-1$. From the equation (\ref{1:eqpa}) we find that
\begin{eqnarray*} \frac{e(q-1)}{h}\left(\frac{q^{m-1}-1}{q-1}\right) = \frac{q^m-1}{q-1}. \end{eqnarray*}
Since
\[\gcd\left(\frac{q^{m-1}-1}{q-1}, \frac{q^m-1}{q-1}\right)=\frac{q^{\gcd(m-1,m)}-1}{q-1}=1, \]
we obtain
\[q^{m-1}-1=q-1 \Longrightarrow m=2. \]
Hence
\[\frac{e(q-1)}{h}= q+1. \]
This implies that $e|(q-1,q+1)=(q-1,2)$. Since $e>1$, we have $e=2$. Now from the above equation we find
\[q=\frac{2+h}{2-h}.\]
The only valid values are $h=1$ and hence $q=3$. So the minimal polynomials of $g^{-1}$ and $(\beta g)^{-1}$ are distinct over $\gf(q)$ except when $q=3,h=1,e=m=2$. $\quad \square$

Next we recall the general result we obtained in \cite[Section 2]{X}.

Denote by $C^{(N,r)}$ the subgroup of $\gf(r)^*$ generated by $\alpha^N$. Since $N|(m,q-1)$, the integer $(r-1)/(q-1)=q^{m-1}+q^{m-2}+\cdots+q+1$ is divisible by $N$, hence $\beta \in C^{(N,r)}$. It is also easy to see that $\gf(q)^* \subset C^{(N,r)}$.

For any $u \in \gf(r)$, define
\begin{eqnarray} \label{2:eta} \eta_u^{(N,r)}=\sum_{z \in C^{(N,r)}}\psi(zu),\end{eqnarray}
where $\psi$ is the canonical additive character of $\gf(r)$, which is given by $\psi(x)=\exp\left(\frac{2 \pi i}{p} \tr_p(x)\right)$, here $\tr_p$ is the trace function from $\gf(r)$ to $\gf(p)$. Obviously $\eta_0^{(N,r)}=\frac{r-1}{N}$. If $u \ne 0$, the term $\eta_u^{(N,r)}$ is called a ``Gaussian period''. Note that the Gaussian periods $\eta_u^{(N,r)}$, $u \ne 0$ depend only on the particular coset of $\gf(r)^*$ with respect to $C^{(N,r)}$ that $u$ belongs to, so there are $N$ such Gaussian periods.

Recall from \cite[Lemma 5]{DL2} (see also \cite{DL1,WT}) that for any $(a,b) \in \gf(r)^2$, the Hamming weight of the codeword $\cc_{(a,b)}$ is given by
\begin{eqnarray} \label{2:weight} \omega\left(\cc_{(a,b)}\right)=\frac{h(r-1)}{q}-\lambda(a,b), \end{eqnarray}
where the ``modified weight'' $\lambda(a,b)$ is defined by
\[\lambda(a,b)=\frac{hN}{eq}\sum_{i=1}^e \eta_{(a+\beta^i b)g^i}^{(N,r)}. \]
It suffices to study $\lambda(a,b)$ only. We have proved in \cite[Section 2]{X} that $\lambda(a,b)$ can attain the following values.

\noindent {\bf Case 1.} $\prod_{i=1}^e\left(a+\beta^i b\right) \ne 0$: For any $c_1,\ldots,c_e \in \gf(r)^*$, we write $\vec{c}=(c_1,\ldots,c_e)$ and define
\[\F(\vec{c})=\left\{(a,b) \in \gf(r)^2: \begin{array}{ll}
\left(a+\beta^i b\right)g^i c_i \in C^{(N,r)} \,\, \forall i \end{array}
\right\}.\]
Then
\begin{eqnarray} \label{2:lambda}
\lambda(a,b)=\frac{hN}{eq}\sum_{i=1}^{e} \eta_{c_i^{-1}}^{(N,r)}, \quad (a,b) \in \F(\vec{c})
\end{eqnarray}
\begin{eqnarray} \label{2:fc}
f(\vec{c}):=\#\F(\vec{c})=\frac{r-1}{N^e}\sum_{\substack{\chi_i^N=\epsilon\\
\chi_1, \ldots, \chi_{e-1}}} f_{\chi_1,\ldots,\chi_{e-1}}(\vec{c}),
\end{eqnarray}
where the sum is over all multiplicative characters $\chi_i$'s of $\gf(r)^*$ such that $\chi_i^N=\epsilon$, $\epsilon$ being the principal character, and
\[f_{\chi_1,\ldots,\chi_{e-1}}(\vec{c})=\prod_{i=1}^{e-1} \chi_i\left(g^{i} (1-\beta^i) c_ic_e^{-1}\right) \sum_{b \in \gf(r)}\prod_{i=1}^{e-1} \chi_i\left(b+\gamma_i\right),\]
\begin{eqnarray} \label{2:gamma}
\gamma_i=\frac{\beta^i}{1-\beta^i}, \quad i=1,2, \ldots, e-1.
\end{eqnarray}

\noindent {\bf Case 2.} $(a,b) \ne (0,0)$ and $a+\beta^t b=0$ for some $t$, $1 \le t \le e$: Then
\begin{eqnarray} \label{2:z2} \lambda\left(-\beta^t b,b\right)=\frac{hN}{eq} \left\{\frac{r-1}{N}+\sum_{\substack{i=1\\ i \ne t}}^e \eta_{bg^i(\beta^i-\beta^t)}^{(N,r)}\right\}, \quad 1 \le t \le e. \end{eqnarray}

\section{The case $e=N=3$: the general setting}

When $e=N=3$, the parameters are
\[\beta=\alpha^{(r-1)/3}, \quad g=\alpha^{(q-1)/h}, \quad 3=\gcd\left(m, \frac{3(q-1)}{h}\right),\quad 3|h \mbox{ and } h|(q-1). \]
Hence $\beta^3=1, 1+\beta+\beta^2=0$. Note that $\beta$ and any $a \in \gf(q)^*$ are both cubic powers in $\gf(r)$. The exact values of the three Gaussian periods $\eta_{u}^{(e,r)}, u=1,\alpha,\alpha^2$ are known (\cite{MY}), which we assume now.

\subsection{Evaluation of $f(\vec{c})$}
We fix a non-trivial cubic character $\rho$ of $\gf(r)^*$. Then all the cubic characters are $\rho,\rho^2$ and $\epsilon=\rho^3$. For $\vec{c}=(c_1,c_2,c_3)$ where $c_1,c_2,c_3 \in \gf(r)^*$, from (\ref{2:lambda}) and (\ref{2:fc}) we have
\begin{eqnarray} \label{3:lambda}
\lambda(a,b)=\frac{h}{q}\sum_{i=1}^{3} \eta_{c_i^{-1}}^{(3,r)}, \quad (a,b) \in \F(\vec{c}),
\end{eqnarray}
\begin{eqnarray} \label{3:fc}
f(\vec{c})&:=&\#\F(\vec{c})=\frac{r-1}{3^3} \sum_{1 \le e_1,e_2 \le 3} \sum_{b \in \gf(r)}\rho\bigl(f_1(b)^{e_1}f_2(b)^{e_2}\bigr),
\end{eqnarray}
where we define
\[f_1(b)=g (1-\beta) c_1c_3^{-1}(b+\gamma_1), \quad f_2(b)=g^2 (1-\beta^2) c_2c_3^{-1}(b+\gamma_2). \]
Here $\gamma_i=\frac{\beta^i}{1-\beta^i}, 1 \le i \le 2$. We see that
\[\gamma_1-\gamma_2=\frac{\beta}{1-\beta^2}. \]
For each $e_1,e_2$, let $\vec{e}=(e_1,e_2)$. Define
\[\A=\{(e_1,e_2): 1 \le e_1,e_2 \le 3\},\quad \A_0=\{(3,3)\},\]
\[\A_1=\{(3,1),(3,2),(3,3)\}, \quad \A_2=\{(1,3),(2,3),(3,3)\},\]
\[\A_3=\{(1,2),(2,1),(3,3)\}, \quad \A_4=\{(1,1),(2,2),(3,3)\},\]
and
\[h_i(\vec{c}):=\sum_{\vec{e} \in \A_i} \sum_{b \in \gf(r)}\rho\bigl(f_1(b)^{e_1}f_2(b)^{e_2}\bigr), \quad 0 \le i \le 4. \]
We can see that
\begin{eqnarray} \label{4:fc} f(\vec{c})=\frac{r-1}{3^3} \Bigl(h_1(\vec{c})+h_2(\vec{c})+h_3(\vec{c})+h_4(\vec{c})-3h_0(\vec{c})\Bigr). \end{eqnarray}
To compute $f(\vec{c})$, it suffices to compute $h_i(\vec{c})$ for each $i$.
Since $\gamma_1 \ne \gamma_2$, we find
\[h_0(\vec{c})=\sum_{\substack{b \in \gf(r)\\
b+\gamma_i \ne 0}}1=r-2. \]
As to $h_1(\vec{c})$, we have
\begin{eqnarray*}
h_1(\vec{c}) &=& \sum_{\vec{e} \in \A_1} \sum_{\substack{b \in \gf(r)\\
b+\gamma_1 \ne 0}} \rho^{e_2}\bigl(f_2(b)\bigr) = \sum_{b \in \gf(r)} \sum_{e=1}^3 \rho^e\bigl(f_2(b)\bigr)-\sum_{e=1}^3 \rho^e\bigl(f_2(-\gamma_1)\bigr).
\end{eqnarray*}
For any $x \in \gf(r)$, define $\delta_3(x)=1$ if $x \in \gf(r)^*$ is a cubic power, and $\delta_3(x)=0$ if otherwise. Then by the orthogonal property of characters we have
\[\delta_3(x)=\frac{1}{3}\sum_{e=1}^3\rho^e(x), \quad x \in \gf(r). \]
Using this we obtain
\[h_1(\vec{c})=\#\{(y,b):y^3=f_2(b), y \ne 0\}-3 \delta_3\bigl(f_2(-\gamma_1)\bigr). \]
From this it is easy to check that
\[h_1(\vec{c})=r-1-3\delta_3\bigl(g^2c_2c_3^{-1}\bigr).\]
Then $h_2(\vec{c})$ is computed in a similar way. We obtain
\[h_2(\vec{c})=r-1-3\delta_3\bigl(gc_1c_3^{-1}\bigr).\]
As for $h_3(\vec{c})$, we obtain
\[h_3(\vec{c})=\sum_{e=1}^3 \sum_{b \in \gf(r)} \rho^e \bigl(f_1(b)f_2(b)^2\bigr)=\sum_{e=1}^3 \sum_{\substack{b \in \gf(r)\\
(b+\gamma_1)(b+\gamma_2) \ne 0}} \rho^e \left(\frac{f_1(b)}{f_2(b)}\right).\]
Make changes of variables we find
\[h_3(\vec{c})=\sum_{e=1}^3 \sum_{\substack{b \in \gf(r)\\
1+(\gamma_1-\gamma_2)b \ne 0\\ b \ne 0}} \rho^e \left(\frac{c_1}{g(1+\beta)c_2}\left(1+(\gamma_1-\gamma_2)b\right)\right).\]
Related to solving the equation
\[y^3=\frac{c_1}{g(1+\beta)c_2}\left(1+(\gamma_1-\gamma_2)b\right), \quad (y,b) \in \gf(r),\]
we find that
\[h_3(\vec{c})=r-1-3 \delta_3 \bigl(g^2c_1c_2^{-1}\bigr). \]
Finally we need to compute $h_4(\vec{c})$. We can write
\[h_4(\vec{c})=\sum_{e=1}^3 \sum_{b \in \gf(r)} \rho \bigl(c_1c_2c_3(b+\gamma_1)(b+\gamma_2)\bigr)=\sum_{e=1}^3 \sum_{b \in \gf(r)} \rho \bigl(c_1c_2c_3b(b+\gamma_1-\gamma_2)\bigr). \]
Let $A$ be the number of solutions $(y,b) \in \gf(r)^2$ such that
\begin{eqnarray} \label{3:c0} C: y^3=c_1c_2c_3b(b+\gamma_1-\gamma_2).\end{eqnarray}
Then $h_4(\vec{c})=A-2$.

\subsubsection{Case 1: $p$ odd.} If $p$ is odd, we can make a change of variables to complete the squares on the right side of (\ref{3:c0}), so that the curve $C$ is transformed into the elliptic curve
\begin{eqnarray} \label{3:codd} E: y^2=x^3-{3(c_1c_2c_3)^4}\,. \end{eqnarray}
The number of $\gf(r)$-points on $E$ can be computed explicitly by using standard tools such as the Gauss sums and the Jacobi sums. For example, following the argument in \cite[Theorem 4, p. 305]{IR} (see also \cite[Exercise 21, p. 63]{KO}), we find that
\[A=r+\chi(-3)\rho\left(c_1c_2c_3\right)J(\rho,\rho)+
\chi(-3)\overline{\rho\left(c_1c_2c_3\right)}\,\overline{J(\rho,\rho)},\]
where $\chi$ is the non-trivial quadratic character of $\gf(r)^*$ and $J(\rho,\rho)$ is the Jacobi sum associated with the character $\rho$. The value $J(\rho,\rho)$ can be evaluated by the Hasse-Davenport relation and \cite[Proposition 8.3.4]{IR}. The explicit result is a little complicated to state, because the choice of one of the two cubic characters $\rho$ has a subtle influence on the exact value of $J(\rho,\rho)$, and $\rho(\alpha)$ could also be any of the two primitive cubic roots of the unity, depending on how we choose the generator $\alpha$ of $\gf(r)^*$. However, the theory about it is well known, so we record the result as follows. Interested readers may refer to \cite[Chapters 8,9,10]{IR} for details.

\begin{lemma} \label{5:le1} Let the assumptions be as before, $r=q^m,q=p^s$, $p$ odd. Let
\[\omega:=\frac{-1+\sqrt{-3}}{2}. \]
Choosing $\rho$ appropriately we have
\[J(\rho,\rho)=(-1)^{ms+1}\pi^{ms},\]
where
\begin{itemize}
\item[1).] If $p \equiv 1 \pmod{3}$, then $\pi \in \Z[\omega]$ is an algebraic integer such that $\pi \bar{\pi}=p$ and $\pi \equiv -1 \pmod{3}$. Identifying $\Z[\omega]/\pi \Z[\omega]$ with the finite field $\FF_p \subset \FF_r$, then $\rho$ is the cubic character of $\FF_r^*$ arising from $\left(\frac{\cdot}{\pi}\right)_3$, the standard cubic residue symbol in $\Z[\omega]/\pi \Z[\omega]$.

\item[2).] If $p \equiv 2 \pmod{3}$, then $s$ is even and $\pi=\sqrt{-p}$. Identifying $\Z[\omega]/p \Z[\omega]$ with the finite field $\FF_{p^2} \subset \FF_r$, then $\rho$ is the cubic character of $\FF_r^*$ arising from the standard cubic residue symbol in $\Z[\omega]/p \Z[\omega]$.
\end{itemize}
\end{lemma}

By Lemma \ref{5:le1} we obtain
\[A=r-\chi(-3)(-1)^{ms}\left(\rho\left(c_1c_2c_3\right)\pi^{ms}+
\rho^2\left(c_1c_2c_3\right)\bar{\pi}^{ms}\right). \]
Therefore
\[h_4(\vec{c})=r-2-\chi(-3)(-1)^{ms}\left(\rho\left(c_1c_2c_3\right)\pi^{ms}+
\rho^2\left(c_1c_2c_3\right)\bar{\pi}^{ms}\right). \]

In summary we obtain
\begin{lemma} \label{4:le3}
If $p$ is odd, then for any $\vec{c}=(c_1,\ldots,c_3)$ where $c_1,\ldots,c_3 \in \gf(r)^*$, we have
\begin{eqnarray*}
f(\vec{c})&=&\frac{r-1}{3^3} \biggl( r+1-3\delta_3(g^2c_2c_3^{-1})-3\delta_3(gc_1c_3^{-1})-3\delta_3(g^2c_1c_2^{-1})\bigr.\\
&&\biggl.-\chi(-3)(-1)^{ms}\bigl(\rho\left(c_1c_2c_3\right)\pi^{ms}+
\rho^2\left(c_1c_2c_3\right)\bar{\pi}^{ms}\bigr)\biggr).
\end{eqnarray*}
\end{lemma}

\subsubsection{Case 2: $p=2$}

On this case, the characteristic of $\gf(r)$ is 2. We find that
\[\gamma_1-\gamma_2=\frac{\beta}{1-\beta}-\frac{\beta^2}{1-\beta^2}=1,\]
so the curve $C$ given in (\ref{3:c0}) is
\begin{eqnarray} \label{6:c2}
C: y^3=c_1c_2c_3b(b+1).
\end{eqnarray}

As in Section 2, $\psi$ is the canonical additive character of $\gf(r)$ given by $\psi(x)=\exp\left( \pi i \tr_2(x)\right)$, here $\tr_2$ is the trace function from $\gf(r)$ to $\gf(2)$. It is easy to see that $x =b^2+b$ for some $b \in \gf(r)$ if and only if $\psi(x)=1$, hence for any $a \in \gf(r)$, the number of solutions for $b \in \gf(r)$ such that $b^2+b=a$ is $1+\psi(a)$. So we find that the number of solutions $(y,b) \in \gf(r)^2$ on the curve $C$ in (\ref{6:c2}) is given by
\[A=\sum_{y \in \gf(r)}1+\psi\left(\frac{y^3}{c_1c_2c_3}\right).\]
From this we obtain
\[A=r+1+3\eta_{(c_1c_2c_3)^2}^{(3,r)}.\]
Therefore
\[h_4(\vec{c})=A-2=r-1+3\eta_{(c_1c_2c_3)^2}^{(3,r)}.\]

In summary we obtain
\begin{lemma} \label{4:le4}
If $p=2$, then for any $\vec{c}=(c_1,\ldots,c_3)$ where $c_1,\ldots,c_3 \in \gf(r)^*$, we have
\begin{eqnarray*}
f(\vec{c})&=&\frac{r-1}{3^3} \biggl( r+2+3\eta_{(c_1c_2c_3)^2}^{(3,r)}-3\delta_3(g^2c_2c_3^{-1})-3
\delta_3(gc_1c_3^{-1})-3\delta_3(g^2c_1c_2^{-1})\biggr). \end{eqnarray*}
\end{lemma}

\subsection{The remaining case}
Next we need to evaluate $\lambda\left(-\beta^tb,b\right)$ in (\ref{2:z2}). For simplicity we adopt a notation: $\lambda \equiv \mu \pmod{\triangle}$ means that $\lambda,\mu \in \gf(r)^*$ and $\frac{\lambda}{\mu}$ is a cube in $\gf(r)^*$.

For $t=1$, it is easy to see that
\[\beta^2-\beta \equiv \beta^3-\beta \equiv \beta-1 \pmod{\triangle}, \]
so we obtain
\begin{eqnarray} \label{4:r1}\lambda\left(-\beta b,b\right)=\frac{h}{q} \left\{\frac{r-1}{3}+ \eta_{bg^2(\beta-1)}^{(3,r)}+\eta_{b(\beta-1)}^{(3,r)}\right\}.\end{eqnarray}

Similarly we find for $t=2$
\begin{eqnarray} \label{4:r2} \lambda\left(-\beta^2 b,b\right)=\frac{h}{q} \left\{\frac{r-1}{3}+ \eta_{bg(\beta-1)}^{(3,r)}+\eta_{b(\beta-1)}^{(3,r)}\right\},\end{eqnarray}
and for $t=3$
\begin{eqnarray} \label{4:r3}\lambda\left(-\beta^3 b,b\right)=\frac{h}{q} \left\{\frac{r-1}{3}+ \eta_{bg(\beta-1)}^{(3,r)}+\eta_{bg^2(\beta-1)}^{(3,r)}\right\}.\end{eqnarray}

\section{The case $e=N=3$: $p \equiv 1 \pmod{3}$}

It is known from (\ref{3:lambda}) that the weight $\lambda(a,b)$ is a simple linear combination of the Gaussian periods. For $u \ne 0$, the Gaussian periods are
\begin{eqnarray} \label{3:g} \eta_u^{(3,r)}=\left\{\begin{array}{lll}
\eta_1^{(3,r)} & \mbox{ if } u \equiv 1 &\pmod{\triangle}, \\
\eta_\alpha^{(3,r)} & \mbox{ if } u \equiv \alpha &\pmod{\triangle}, \\
\eta_{\alpha^2}^{(3,r)} & \mbox{ if } u \equiv \alpha^2 &\pmod{\triangle}.
\end{array} \right.\end{eqnarray}
The exact values of the three Gaussian periods are known (\cite{MY}), which we assume now and will make explicit later. Basically, when $p \equiv 1 \pmod{3}$, the three Gaussian periods are all distinct.

We summarize the argument as follows: each $c_i (\ne 0)$, $i=1,2,3$ has three distinct values: $c_i \equiv 1,\alpha, \alpha^2 \pmod{\triangle}$, which result in three Gaussian periods $\eta_{c_i^{-1}}^{(3,r)}$, so $\lambda(a,b)$ from the equation (\ref{3:lambda}) has at most 10 different values, and the value depends only on the vector $\vec{c}=(c_1,c_2,c_3)$ for which $(a,b) \in \F(\vec{c})$. Moreover, for each such vector $\vec{c}$, the number $\#\F(\vec{c})=f(\vec{c})$ is given by Lemma \ref{4:le3}. We also notice that when $p \equiv 1 \pmod{3}$, then $\left(\frac{-3}{p}\right)=1$, here $\left(\frac{\cdot}{p}\right)$ is the Legendre symbol with respect to $p$, hence
\[\chi(-3)=\left(\frac{-3}{p}\right)^{ms}=1. \]
The final results will depend on whether or not $3|\frac{q-1}{h}$.

\subsection{The Proof of Table \ref{1:t3}}
Assume first that $3|\frac{q-1}{h}$, then $g=\alpha^{(q-1)/h}$ is a cube in $\gf(r)$. By Lemma \ref{4:le3} we have
\begin{eqnarray*}
f(\vec{c})&=&\frac{r-1}{3^3} \biggl( r+1-3\delta_3(c_2c_3^{-1})-3\delta_3(c_1c_3^{-1})-3\delta_3(c_1c_2^{-1})\bigr.\\
&&\biggl.-(-1)^{ms}\bigl(\rho\left(c_1c_2c_3\right)\pi^{ms}+
\rho^2\left(c_1c_2c_3\right)\bar{\pi}^{ms}\bigr).
\end{eqnarray*}

For any $u=1, \alpha$ or $\alpha^2$, if $c_1 \equiv c_2 \equiv c_3 \equiv u \pmod{\triangle}$, then the ``modified weight'' $\lambda(a,b)= \frac{h}{q}3 \eta_{u^{-1}}^{(3,r)}= \frac{3h}{q}\eta_{u^{-1}}^{(3,r)}$, and the total number of such $(a,b)$'s counted in this $\F(\vec{c})$ is given by
\[f(\vec{c})=\frac{r-1}{3^3} \biggl(r-8-(-1)^{ms} \bigl(\pi^{ms}+\bar{\pi}^{ms}\bigr)\biggr).\]

If $c_1,c_2,c_3$ are all distinct modulo $\triangle$, then $\lambda(a,b)= \frac{h}{q} \left(\eta_{1}^{(3,r)}+\eta_{\alpha}^{(3,r)}+
\eta_{\alpha^2}^{(3,r)}\right)=-\frac{h}{q}$, and the number of such $(a,b)$'s counted in this $\F(\vec{c})$ is given by
\[f(\vec{c})=\frac{r-1}{3^3} \biggl(r+1-(-1)^{ms} \bigl(\pi^{ms}+\bar{\pi}^{ms}\bigr)\biggr).\]
The total number of such $\vec{c}$'s is 6.

Let $(i,j,k)$ be a permutation of $(1,2,3)$. If $c_i \equiv c_j \equiv 1 \pmod{\triangle}$ and $c_k \equiv \alpha \pmod{\triangle}$, then $\lambda(a,b)= \frac{h}{q} \left(2\eta_{1}^{(3,r)}+\eta_{\alpha^2}^{(3,r)}\right)$, and the number of such $(a,b)$'s counted is given by
\[f(\vec{c})=\frac{r-1}{3^3} \biggl(r-2-(-1)^{ms} \bigl(\rho(\alpha)\pi^{ms}+\rho^2(\alpha)\bar{\pi}^{ms}\bigr)\biggr).\]
The total number of such $\vec{c}$'s is 3. On the other hand, if $c_i \equiv c_j \equiv 1 \pmod{\triangle}$ and $c_k \equiv \alpha^2 \pmod{\triangle}$, then $\lambda(a,b)= \frac{h}{q} \left(2\eta_{1}^{(3,r)}+\eta_{\alpha}^{(3,r)}\right)$, and the number of such $(a,b)$'s counted is given by
\[f(\vec{c})=\frac{r-1}{3^3} \biggl(r-2-(-1)^{ms} \bigl(\rho^2(\alpha)\pi^{ms}+\rho(\alpha)\bar{\pi}^{ms}\bigr)\biggr).\]
The total number of such $\vec{c}$'s is also 3.

Similarly we consider that cases that
\begin{itemize}
\item $c_i \equiv c_j \equiv \alpha \pmod{\triangle}$ and $c_k \equiv 1$ or $\alpha^2 \pmod{\triangle}$;

\item $c_i \equiv c_j \equiv \alpha^2 \pmod{\triangle}$ and $c_k \equiv 1$ or $\alpha \pmod{\triangle}$;
\end{itemize}
and we can obtain similar results. This yields the first ten weights $\lambda(a,b)$ and the corresponding frequencies in Table \ref{1:t3}.

We also need to count the weight $\lambda(a,b)$ and its frequency coming from $(a,b)$'s such that $a=-\beta^t b$ for some $t$, $1 \le t \le 3$. Such cases are treated in (\ref{4:r1})--(\ref{4:r3}). Since $3|\frac{q-1}{h}$, $g$ is a cube in $\gf(r)$. Considering $t=1$ and $\lambda\left(-\beta b,b\right)$ in (\ref{4:r1}), we find that
\[ \lambda\left(-\beta b,b\right)=\frac{h}{q}\left\{\frac{r-1}{3}+2 \eta_{b(\beta-1)}^{(3,r)}\right\}.\]
This is $\frac{h}{q}\left\{\frac{r-1}{3}+2 \eta_{u}^{(3,r)}\right\}$ with frequency $(r-1)/3$ for any $u=1,\alpha,\alpha^2$ respectively when $b \in \gf(r)^*$ varies. The results for $t=2,3$ are the same, and this yields the weights $\lambda(a,b)$ and frequencies from lines 11--13 in Table \ref{1:t3}. The last line of Table \ref{1:t3} comes from $(a,b)=(0,0)$, which corresponds to the codeword with Hamming weight zero. The completes the proof of Table \ref{1:t3}. $\square$

\subsection{The Proof of Table \ref{1:t4}}

Assume that $3 \nmid \frac{q-1}{h}$, so that $g \equiv \alpha$ or $\alpha^2 \pmod{\triangle}$. The analysis is similar in either case.

For any $u=1, \alpha$ or $\alpha^2$, if $c_1 \equiv c_2 \equiv c_3 \equiv u \pmod{\triangle}$, then the ``modified weight'' $\lambda(a,b)=  \frac{3h}{q}\eta_{u^{-1}}^{(3,r)}$, and the total number of such $(a,b)$'s counted in this $\F(\vec{c})$ is given by
\[f(\vec{c})=\frac{r-1}{3^3} \biggl(r+1-(-1)^{ms} \bigl(\pi^{ms}+\bar{\pi}^{ms}\bigr)\biggr).\]

If $c_1,c_2,c_3$ are all distinct modulo $\triangle$, then $\lambda(a,b)= -\frac{h}{q}$, and the number of such $(a,b)$'s counted in this $\F(\vec{c})$ is given by
\begin{eqnarray*}f(\vec{c})&=&\frac{r-1}{3^3} \biggl(r+1-(-1)^{ms} \bigl(\pi^{ms}+\bar{\pi}^{ms}\bigr)\\
&&\biggl. -3\delta_3\left(g^2c_2c_3^{-1}\right)-3\delta_3\left(gc_1c_3^{-1}\right)-3\delta_3\left(g^2c_1c_2^{-1}\right)\biggr).\end{eqnarray*}
The total number of such $\vec{c}$'s is 6. Collecting all such $\vec{c}$'s into the set $\Ta$, it is easy to check that
\[\sum_{\vec{c} \in \Ta}\delta_3\left(g^2c_2c_3^{-1}\right)=\sum_{\vec{c} \in \Ta}\delta_3\left(gc_1c_3^{-1}\right)=\sum_{\vec{c} \in \Ta}\delta_3\left(g^2c_1c_2^{-1}\right)=3.\]
Hence the total number of such $(a,b)$'s inside one of those $\F(\vec{c})$'s is given by
\begin{eqnarray*}\sum f(\vec{c})&=&\frac{r-1}{3^2} \biggl(2r-7-2(-1)^{ms} \bigl(\pi^{ms}+\bar{\pi}^{ms}\bigr)\biggr).\end{eqnarray*}

Let $(i,j,k)$ be a permutation of $(1,2,3)$. If $c_i \equiv c_j \equiv 1 \pmod{\triangle}$ and $c_k \equiv \alpha \pmod{\triangle}$, then $\lambda(a,b)= \frac{h}{q} \left(2\eta_{1}^{(3,r)}+\eta_{\alpha^2}^{(3,r)}\right)$, and the number of such $(a,b)$'s counted is given by
\begin{eqnarray*}f(\vec{c})&=&\frac{r-1}{3^3} \biggl(r+1-\chi(-3)(-1)^{ms} \bigl(\rho(\alpha)\pi^{ms}+\rho^2(\alpha)\bar{\pi}^{ms}\bigr)\\
&&\biggl. -3\delta_3\left(g^2c_2c_3^{-1}\right)-3\delta_3\left(gc_1c_3^{-1}\right)-3\delta_3\left(g^2c_1c_2^{-1}\right)\biggr).\end{eqnarray*}
The total number of such $\vec{c}$'s is 3. Collecting all such $\vec{c}$'s into the set $\Tb$, it is easy to check that
\[\sum_{\vec{c} \in \Tb}\delta_3\left(g^2c_2c_3^{-1}\right)=\sum_{\vec{c} \in \Tb}\delta_3\left(gc_1c_3^{-1}\right)=\sum_{\vec{c} \in \Tb}\delta_3\left(g^2c_1c_2^{-1}\right)=1.\]
Hence the total number of such $(a,b)$'s inside one of those $\F(\vec{c})$'s is given by
\begin{eqnarray*}\sum f(\vec{c})&=&\frac{r-1}{3^2} \biggl(r-2-(-1)^{ms} \bigl(\rho(\alpha)\pi^{ms}+\rho^2(\alpha)\bar{\pi}^{ms}\bigr)\biggr).\end{eqnarray*}
For all the other cases the results are similar. This yields the first 10 weights $\lambda(a,b)$ and the corresponding frequencies in Table \ref{1:t4}.

We also need to count the weight $\lambda(a,b)$ and its frequency coming from $(a,b)$'s such that $a=-\beta^t b$ for some $t$, $1 \le t \le 3$. Such cases are treated in (\ref{4:r1})--(\ref{4:r3}). Since $3 \nmid \frac{q-1}{h}$, $g$ is not a cube in $\gf(r)$. Considering $t=1$ and $\lambda\left(-\beta b,b\right)$ in (\ref{4:r1}), we find that
\[ \lambda\left(-\beta b,b\right)=\frac{h}{q}\left\{\frac{r-1}{3}+ \eta_{bg^2(\beta-1)}^{(3,r)}+\eta_{b(\beta-1)}^{(3,r)}\right\}.\]
The right hand side could be $\frac{h}{q}\left\{\frac{r-1}{3}+ \eta_{1}^{(3,r)}+\eta_{\alpha}^{(3,r)}\right\}$, $\frac{h}{q}\left\{\frac{r-1}{3}+ \eta_{1}^{(3,r)}+\eta_{\alpha^2}^{(3,r)}\right\}$ or $\frac{h}{q}\left\{\frac{r-1}{3}+ \eta_{\alpha}^{(3,r)}+\eta_{\alpha^2}^{(3,r)}\right\}$, each of which appears with frequency $(r-1)/3$ when $b \in \gf(r)^*$ varies. The results for $t=2,3$ are the same. This yields the weights $\lambda(a,b)$ and frequencies from lines 11--13 in Table \ref{1:t4}. The last line of Table \ref{1:t4} comes from $(a,b)=(0,0)$, which corresponds to the codeword with Hamming weight zero. This completes the proof of Table \ref{1:t4}. $\square$

\subsection{The Gaussian Periods}

Finally we shall determine explicitly the three Gaussian periods $\eta_{1}^{(3,r)},\eta_{\alpha}^{(3,r)}$ and $\eta_{\alpha^2}^{(3,r)}$. This could be provided by \cite[Theorem 22]{MY}, however, it seems the results do not distinguish the value $\eta_{\alpha}^{(3,r)}$ from $\eta_{\alpha^2}^{(3,r)}$, so we use the basic property \cite[Proposition 1]{MY} to compute the values by ourselves. It works in general situation, but for simplicity, we stick to our notation and the special case that $e=N=3$ and $p \equiv 1 \pmod{3}$.

First, for $k=0,1,2$, define
\[G_k:=\sum_{x \in \gf(r)} \psi(\alpha^k x^3), \]
here as in Section 2, $\psi$ is the canonical additive character of $\gf(r)$. It is easy to see that
\[G_k=3\eta_{\alpha^k}^{(3,r)}+1, \quad k=0,1,2. \]
So it is enough to find the values $G_k$. \cite[(g) of Proposition 1]{MY} states the relation
\[G_k=\sum_{j=1}^2\rho(\alpha)^{-jk} \tau(\rho^j), \]
where $\rho$ is the cubic character of $\gf(r)$ arising from $\chi_{\pi}(\cdot):=\left(\frac{\cdot}{\pi}\right)_3$ that we have chosen from Lemma 1, and for any  multiplicative character $\zeta$ of $\gf(r)$, $\tau(\zeta)$ is the Gauss sum given by
\[\tau(\zeta):=\sum_{x \in \gf(r)^*} \zeta(x) \psi(x). \]
Because of the careful choice of $\pi$, we know that (\cite[Section 4, Chapter 9]{IR})
\[J(\chi_{\pi},\chi_{\pi})=\pi, \quad \tau(\chi_{\pi})^3=p\pi, \]
where $J(\chi_{\pi},\chi_{\pi})$ and $\tau(\chi_{\pi})$ are the Jacobi sum and Gauss sum defined on $\Z[\omega]/\pi\Z[\omega]$, which is identified naturally as $\gf(p)$. Since $r=p^{sm}$ and $3|sm$, by the Davenport-Hasse relation (see \cite{IR} or \cite[Proposition 18]{MY}), we find that
\[\tau(\rho)=(-1)^{ms+1}(\tau(\chi_{\pi}))^{ms}=(-1)^{ms+1}
\left(\tau(\chi_{\pi})^3\right)^{ms/3}=(-1)^{ms+1}r^{1/3}\pi^{sm/3}.\]
We also derive
\[\tau(\rho^2)=\overline{\tau(\rho)}=
(-1)^{ms+1}r^{1/3}\bar{\pi}^{sm/3}.\]
The values $\tau(\rho)$ and $\tau(\rho^2)$ can be used to evaluate $G_k$'s, which in turn provide explicit evaluations of $\eta_{\alpha^k}^{(3,r)}$'s as claimed in Theorem \ref{thm2}. Now the proof of Theorem \ref{thm2} is complete. $\square$

\section{The case $e=N=3$: $p \equiv 2 \pmod{3}$}

\subsection{}

We first argue that the cases $p=2$ and $p \equiv 2 \pmod{3}$, $p$ odd, can be brought together.

First, when $p \equiv 2 \pmod{3}$, whether $p$ is odd or even, then $2|s$. The values of the Gaussian periods are described as follows.
\begin{lemma}[Proposition 20, \cite{MY}] \label{gaussp1}
Assume that $p \equiv 2 \pmod{3}$, $r=p^{sm}$, $6|sm$. Define \[G_u=3\eta_u^{(3,r)}+1, \quad u \in \gf(r)^*. \]
Then
\[G_1=-2(-1)^{ms/2}\sqrt{r}, \quad G_{\alpha}=G_{\alpha^2}=
(-1)^{ms/2}\sqrt{r}.\]
\end{lemma}
So in either case, there are two distinct values in the Gaussian periods, and the formulas are the same.

Second, when $p \equiv 2 \pmod{3}$ and $p$ is odd, then $\pi=\sqrt{-p}$. Since \[\frac{r-1}{p-1}=p^{sm-1}+p^{sm-2}+\ldots+p+1 \equiv 0 \pmod{2}, \]
and $\FF_p^*$ is generated by $\alpha^{(r-1)/(p-1)}$, any $a \in \FF_p^*$ is a perfect square in $\gf(r)$, hence $\chi(-3)=1$. For any $\vec{c}=(c_1,c_2,c_3)$ where $c_1,c_2,c_3 \in \gf(r)^*$, the formulas for $f(\vec{c})$ given by Lemma \ref{4:le3} can be simplified as
\begin{eqnarray} \label{5:fcp}
f(\vec{c})&=&\left\{\begin{array}{c}\frac{r-1}{3^3} \biggl( r+1-3\delta_3(g^2c_2c_3^{-1})-3\delta_3(gc_1c_3^{-1})-3\delta_3(g^2c_1c_2^{-1})\bigr.\\
\biggl.-(-1)^{ms/2}\sqrt{r}\bigl(\rho\left(c_1c_2c_3\right)+
\rho^2\left(c_1c_2c_3\right)\bigr)\biggr).\end{array} \right.
\end{eqnarray}
Using
\[\rho\left(c_1c_2c_3\right)+
\rho^2\left(c_1c_2c_3\right)=\left\{\begin{array}{cc}
2& \mbox{ if } c_1c_2c_3 \in C^{(3,r)},\\
-1& \mbox{ if } c_1c_2c_3 \not \in C^{(3,r)},
\end{array}\right.\]
the formulas for $f(\vec{c})$ can be simplified further.

On the other hand, if $p=2$, using the Gaussian periods described by Lemma \ref{gaussp1}, it is easy to find that for any $\vec{c}=(c_1,c_2,c_3)$ where $c_1,c_2,c_3 \in \gf(r)^*$, the formula for $f(\vec{c})$ given by Lemma \ref{4:le4} is the same as the one for $f(\vec{c})$ given in (\ref{5:fcp}). Therefore the results for $p=2$ can be included into the case $p \equiv 2 \pmod{3}$.

\subsection{} The argument for $p \equiv 2 \pmod{3}$, $p$ odd, is exactly the same as the previous section for $p \equiv 1 \pmod{3}$, and the results are summarized in Tables \ref{1:t3}-\ref{1:t4}, except that we have to take into account of the extra conditions
\[\eta_{\alpha}^{(3,r)}=\eta_{\alpha^2}^{(3,r)}, \quad \eta_{1}^{(3,r)}+2\eta_{\alpha}^{(3,r)}=-1.\]
In other words, some of the different weights in Table \ref{1:t3}-\ref{1:t4} turn out to be the same if $p \equiv 2 \pmod{3}$. Taking care of those repeated weights, we have completed the proof of Theorem \ref{thm1}.

\end{document}